\documentclass[11pt,twoside]{article}
\usepackage[plainpages=false]{hyperref}
\usepackage{amsfonts,latexsym,rawfonts,amsmath,amssymb}
\usepackage{amsmath,amssymb,amsfonts,latexsym,lscape,rawfonts}
\newenvironment{proof}[1][Proof]{\begin{trivlist}
\item[\hskip \labelsep {\bfseries #1}]}{\end{trivlist}}

\newcommand{\twoten}[3]{#1_{#2\bar #3}}

\newcommand{\four}[5]{#1_{#2\bar #3#4\bar #5}}

\newcommand{\pd}[2]{\frac {\partial #1}{\partial #2}}

\newcommand{\al}{\alpha}
\newcommand{\bb}{\beta}
\newcommand{\la}{\lambda}
\newcommand{\La}{\Lambda}

\newcommand{\Na}{\nabla}

\newcommand{\ee}{\epsilon}

\newcommand{\beqs}{\begin{eqnarray*}}
\newcommand{\eeqs}{\end{eqnarray*}}
\newcommand{\beqn}{\begin{eqnarray}}
\newcommand{\eeqn}{\end{eqnarray}}
\newcommand{\beqa}{\begin{array}}
\newcommand{\eeqa}{\end{array}}

\def\lra{\longrightarrow}
\def\La{\Lambda}
\def\la{\lambda}

\def\p{\prime}

\def\bc{\begin{center}}
\def\ec{\end{center}}

\def\p{\partial}
\def\b{\bar}

\def\S{\Sigma}

\def\C{{\mathbb C}}

\def\B{{\mathfrak B}}
\def\C{{\mathfrak C}}

\def\M{{\mathfrak M}}

\def\cR{{\mathcal R}}
\def\cS{{\mathcal S}}

\def\RR{{\mathbb R}}

\newtheorem{prop}{Proposition}[section]
\newtheorem{theo}[prop]{Theorem}
\newtheorem{lem}[prop]{Lemma}
\newtheorem{claim}[prop]{Claim}
\newtheorem{cor}[prop]{Corollary}
\newtheorem{rem}[prop]{Remark}

\newtheorem{defi}[prop]{Definition}

\def\begeq{\begin{equation}}
\def\endeq{\end{equation}}
\def\and{\quad{\rm and}\quad}

\let\lra=\longrightarrow

\def\mapright\#1{\,\smash{\mathop{\lra}\limits^{\#1}}\,}

\title{{\bf\Large{The K\"ahler-Ricci flow on K\"ahler manifolds with 2 traceless
bisectional curvature operator}}}

\author{X. X. Chen, H. Li}
\begin{document}
\bibliographystyle{plain}
\maketitle
\begin{center} {\bf  Dedicated to Professor W. Y. Ding for his 60's birthday}
\end{center}

\tableofcontents
\section{Introduction}
In 1982, in a  famous paper  \cite{[8]}, R. Hamilton  proved
that in a 3 dimensional compact  manifold, if the initial metric has
positive Ricci curvature, then this positivity condition will be
preserved under the Ricci flow.  He showed that the underlying manifold must be diffeomorphic to the standard $S^3\;$ or its finite quotient.  Following this paper,  there are intensive active
researches on Ricci flow, and many works are devoted to study when certain convex cones of
curvature pinching condition is preserved by the Ricci flow. In
\cite{[9]}, R. Hamilton proved that the positive curvature operator
is preserved under the Ricci flow in all dimensions.
 H. Chen\cite{[4]} further showed that, a weaker notion, that the sum of
any two eigenvalues is positive, is preserved under the Ricci flow.
 In the K\"ahler setting, it is well known that the positive bisectional curvature
is preserved under the K\"ahler Ricci flow
 through the work of S. Bando \cite{[1]} for complex dimension $n=3$, and later N. Mok \cite{[12]} for general
 dimension.  Following the argument of N. Mok, in an unpublished work of Cao-Hamilton, they proved that the orthogonal bisectional curvature is preserved
 under the K\"ahler Ricci flow.  There are other convex cones of curvature pinching conditions which
 are  preserved, for instance
\cite{ [2a]} and \cite{ [11a]}.  A more complete reference on this topic
 can be found in \cite{[10]} .\\

Analyzing the evolution equation (2.4) of the Ricci tensor,  it is somewhat unfortunate that
the parabolic Laplacian of the Ricci
tensor involves the full sectional curvature.  It is then not
surprise that we only know the positivity of Ricci tensor is
preserved in real dimension $2$ and $3$ by the earlier work of R.
Hamilton. A counter example to the possible extension of R. Hamilton's result on Ricci tensor in high dimension
seems to be difficult to construct.  More recently, Ni Lei
\cite{[11]} constructed first counter example to this in
Riemanninian setting where the positivity of Ricci curvature is not preserved by the Ricci flow.
Very recently, Dan Knopf \cite{[7a]} constructed similar counter example
in the K\"ahler setting. However, both examples  are in non-compact
manifold.   Therefore, it is still an open question whether or not
positive Ricci curvature is preserved under the Ricci flow in the
case of compact manifold. In particular, in the case of compact
K\"ahler manifold,  there might be some hope
that some form of lower bound of Ricci curvature will be preserved in \cite{[5]} where the first named author showed, along with other results\footnote{In \cite{[5]} , the first named author proved that,  any irreducible K\"ahler manifold where the positive orthogonal bisectional curvature is preserved under the K\"ahler Ricci flow,
must be biholomorphic to ${\mathbb C} {\mathbb P}^n.\;$},  that any metric with positive orthogonal bisectional curvature, even
a negative lower bound of Ricci curvature, is preserved and improved under the K\"ahler Ricci flow.\\

  In a compact K\"ahler manifold $X$, the bisectional
curvature tensor acts as a symmetric bilinear form on  the space of
$(1,1)$ form (which we will denote as $\Lambda^{1,1}(X)).\;$
Furthermore, this action (of traceless part of the bisectional curvature)  preserves the traceless part of this space
(which we will denote as $\Lambda^{1,1}_0(X)$).   In a recent paper
by Phong-Sturm \cite{[15]}, they observed that, the condition
that
 the sum of any two eigenvalues of the traceless bisectional curvature operator
is positive, is preserved under the K\"ahler Ricci flow in complex
dimension 2. Note that this condition is different from the
condition used by H. Chen, even though the main idea of proof is
very similar.  The main theorem they proved in \cite{[15]} is that,
if this curvature condition hold, then positive Ricci curvature will
be preserved under the  K\"ahler Ricci flow in complex surface. The proof there is  complicated,
largely due to the fact that the action of curvature operator in $\Lambda^{1,1}_0(X)$
 is very complicated. \\

  The  2-positive
traceless bisectional curvature is certainly different to the
popular notion of positive bisectional curvature. For instance, when
this curvature condition holds, the Ricci curvature might not be
positive.  In \cite{[6]} \cite{[7]}, the first named author and G.
Tian studied the convergence of K\"ahler Ricci flow in K\"ahler
Einstein manifold where the initial metric has positive bisectional
curvature and showed that the K\"ahler Ricci flow must converges to
the Fubni-Study metric exponentially over the flow.  The present
work can be viewed as a continuation of \cite{[6]} \cite{[7]} in the
sense that the curvature condition is relaxed in some subtle way.
However, this type of curvature condition was indeed studied first by H. Chen in \cite{[4]}.
%More recently, in an elegantly written paper \cite{[15]}, Jacob-Phong studied these curvature
%conditions in the special case of complex K\"ahler surface.
 %The interest of the first named author in this type of special curvature conditions was
%certainly re-invigorated by this  paper \cite{[15]}.\\
Following  \cite{[4]} \cite{[15]} , we
study systematically geometrical properties of this 2-positive
traceless bisectional curvature operator on any K\"ahler manifold.
 Our first result is:

\begin{theo} Let $X$ be a compact K\"ahler manifold of dimension
$n\geq 2.\;$ Along the K\"ahler Ricci flow, we have

 \begin{enumerate}
\item If the initial metric has non-negative traceless bisectional
curvature operator, then the evolved metrics also have
non-negative traceless bisectional curvature operator. If it is
positive at one point initially, then it is positive everywhere
for all $t>0.\;$ \item If the initial metric has 2-non-negative
traceless bisectional curvature operator, then the evolved metrics
also have 2-non-negative traceless bisectional curvature. If it is
positive at one point initially, then it is positive everywhere
for all $t>0.\;$

\end{enumerate}
Under either of these two conditions, the positivity of Ricci
tensor is preserved under the K\"ahler Ricci flow.
\end{theo}

The relation of 2-positive traceless bisectional curvature with
 the notion of positive orthogonal bisectional curvature is much more subtle.  They are defined
 in a completely different manner and the action of bisectional curvature operator on the space of
 (1,1) forms is very complicated.  It is hard to visualize what 2 positive traceless bisectional curvature really is.     A somewhat surprising result
 we prove in this paper is that  (c.f. Theorem 1.2
below)  any K\"ahler metric  which has 2-positive traceless
bisectional curvature, must also have positive orthogonal
bisectional curvature.  The last part of the preceding theorem
follows directly from the application of Hamilton's maximal
principle for tensors to the evolution equation of the Ricci tensor.
Comparing to the main theorem in \cite{[15]}, our theorem is for all
dimensional and our proof  is simpler and more straightforward.

  \begin{theo} In a K\"ahler manifold with 2-non-negative traceless bisectional curvature operator, then
 the orthogonal bisectional curvature must be non-negative.    If the scalar curvature is uniformly bounded,
 then   the bisectional curvature
 is uniformly bounded.  Moreover, if we assume that the traceless bisectional curvature operator is non-negative, then
the sum of any two eigenvalues of the Ricci tensor is
non-negative.
 \end{theo}

 \begin{rem}  The proof that the scalar curvature controls the full bisectional curvature is sophiscated,  and somewhat lengthy.
 In the special case of complex surface, similar estimate is derived also in \cite{[16]}. In an unpublished work of G. Perelman, the scalar curvature is uniformly bounded
 along the K\"ahler Ricci flow.  Combining this with Theorem 1.2, we conclude that the full bisectional curvature is uniformly bounded over
 the K\"ahler Ricci flow when the initial metric has 2 positive traceless bisectional curvature.
 \end{rem}

Following Remark 1.3 and a general theorem on the K\"ahler Ricci flow (c.f. \cite{[17]}, \cite{[18]})
we arrive the following
\begin{cor}  Along the K\"ahler Ricci flow, if the initial metric has 2-positive traceless bisectional curvature operator, then the flow converges by sequence to some K\"ahler Ricci Soliton in
the limit in the sense of Cheeger-Gromov.
 \end{cor}

Similar results are also proved by Jacob -Phong \cite{[16]} in
special case of complex  K\"ahler surfaces.

\vskip1cm
\section{ Basic  K\"ahler geometry}
\subsection{Setup of notations}

Let $X$ be an $n$-dimensional compact K\"ahler manifold. A
K\"ahler metric can be given by its K\"ahler form $\omega$ on $X$.
In local coordinates $z_1, \cdots, z_n$, this $\omega$ is of the
form
\[
\omega = \sqrt{-1} \displaystyle \sum_{i,j=1}^n\;g_{i
\overline{j}} d\,z^i\wedge d\,z^{\overline{j}}  > 0,
\]
where $\{g_{i\overline {j}}\}$ is a positive definite Hermitian matrix function.
The K\"ahler condition requires that $\omega$ is a closed positive
(1,1)-form. In other words, the following holds
\[
 {{\partial g_{i \overline{k}}} \over
{\partial z^{j}}} =  {{\partial g_{j \overline{k}}} \over
{\partial z^{i}}}\qquad {\rm and}\qquad {{\partial g_{k
\overline{i}}} \over {\partial z^{\overline{j}}}} = {{\partial
g_{k \overline{j}}} \over {\partial
z^{\overline{i}}}}\qquad\forall\;i,j,k=1,2,\cdots, n.
\]
The K\"ahler metric corresponding to $\omega$ is given by
\[
 \sqrt{-1} \;\displaystyle \sum_1^n \; {g}_{\alpha \overline{\beta}} \;
d\,z^{\alpha}\;\otimes d\, z^{ \overline{\beta}}.
\]
For simplicity, in the following, we will often denote by $\omega$
the corresponding K\"ahler metric. The K\"ahler class of $\omega$
is its cohomology class $[\omega]$ in $H^2(X,\RR).\;$ By the Hodge
theorem, any other K\"ahler metric in the same K\"ahler class is
of the form
\[
\omega_{\varphi} = \omega + \sqrt{-1} \displaystyle \sum_{i,j=1}^n\;
{{\partial^2 \varphi}\over {\partial z^i \partial z^{\overline{j}}}}
> 0
\]
for some real valued function $\varphi$ on $X.\;$ The functional
space in which we are interested (often referred as the space of
K\"ahler potentials) is
\[
{\cal P}(X,\omega) = \{ \varphi \;\mid\; \omega_{\varphi} = \omega
+ \sqrt{-1}
 {\partial} \overline{\partial} \varphi > 0\;\;{\rm on}\; X\}.
\]
Given a K\"ahler metric $\omega$, its volume form  is
\[
  \omega^n = {1\over {n!}}\;\left(\sqrt{-1} \right)^n \det\left(g_{i \overline{j}}\right)
 d\,z^1 \wedge d\,z^{\overline{1}}\wedge \cdots \wedge d\,z^n \wedge d\,z^{\overline{n}}.
\]
Its Christoffel symbols are given by
\[
  \Gamma^k_{i\,j} = \displaystyle \sum_{l=1}^n\;g^{k\overline{l}} {{\partial g_{i \overline{l}}} \over
{\partial z^{j}}} ~~~{\rm and}~~~ \Gamma^{\overline{k}}_{\overline{i}\,\overline{j}} =
\displaystyle \sum_{l=1}^n\;g^{\overline{k}l} {{\partial g_{l \overline{i}}} \over
{\partial z^{\overline{j}}}}, \qquad\forall\;i,j,k=1,2,\cdots n.
\]
The curvature tensor is
\[
 R_{i \overline{j} k \overline{l}} = - {{\partial^2 g_{i \overline{j}}} \over
{\partial z^{k} \partial z^{\overline{l}}}} + \displaystyle \sum_{p,q=1}^n g^{p\overline{q}}
{{\partial g_{i \overline{q}}} \over
{\partial z^{k}}}  {{\partial g_{p \overline{j}}} \over
{\partial z^{\overline{l}}}}, \qquad\forall\;i,j,k,l=1,2,\cdots n.
\]
We say that $\omega$ is of nonnegative bisectional curvature if
\[
 R_{i \overline{j} k \overline{l}} v^i v^{\overline{j}} w^k w^{\overline{l}}\geq 0
\]
for all non-zero vectors $v$ and $w$ in the holomorphic tangent
bundle of $X$. The bisectional curvature and the curvature tensor
can be mutually determined. The Ricci curvature of $\omega$ is
locally given by
\[
  R_{i \overline{j}} = - {{\partial}^2 \log \det (g_{k \overline{l}}) \over
{\partial z_i \partial \bar z_j }} .
\]
So its Ricci curvature form is
\[
  {\rm Ric}(\omega) = \sqrt{-1} \displaystyle \sum_{i,j=1}^n \;R_{i \overline{j}}(\omega)
d\,z^i\wedge d\,z^{\overline{j}} = -\sqrt{-1} \partial \overline{\partial} \log \;\det (g_{k \overline{l}}).
\]
It is a real, closed (1,1)-form. Recall that $[\omega]$ is called
a canonical K\"ahler class if this Ricci form is cohomologous to
$\lambda \;\omega,\; $ for some constant $\lambda.$ In our
setting, we require $\lambda = 1.\;$

 \subsection{The K\"ahler Ricci flow}

    Now we assume that the first Chern class $c_1(X)$ is positive.
The normalized Ricci flow (c.f. \cite{[8]} and \cite{[9]}) on a
K\"ahler manifold $X$ is of the form
\begin{equation}
  {{\partial g_{i \overline{j}}} \over {\partial t }} = g_{i \overline{j}}
  - R_{i \overline{j}}, \qquad\forall\; i,\; j= 1,2,\cdots ,n.
\label{eq:kahlerricciflow}
\end{equation}
If we choose the initial K\"ahler metric $\omega$ with $c_1(X)$ as
its K\"ahler class. The flow (2.1) preserves the K\"ahler class
$[\omega]$. It follows that on the level of K\"ahler potentials,
the Ricci flow becomes
\begin{equation}
   {{\partial \varphi} \over {\partial t }} =  \log {{\omega_{\varphi}}^n \over {\omega}^n } + \varphi - h_{\omega} ,
\label{eq:flowpotential}
\end{equation}
where $h_{\omega}$ is defined by
\[
  {\rm Ric}(\omega)- \omega = \sqrt{-1} \partial \overline{\partial} h_{\omega}, \; {\rm and}\;\displaystyle \int_X\;
  (e^{h_{\omega}} - 1)  {\omega}^n = 0.
\]
Then the evolution equation for bisectional curvature is

\begin{eqnarray}{{\partial }\over {\partial t}} R_{i \overline{j} k
\overline{l}} & = & \bigtriangleup R_{i \overline{j} k
\overline{l}} + R_{i \overline{j} p \overline{q}} R_{q
\overline{p} k \overline{l}} - R_{i \overline{p} k \overline{q}}
R_{p \overline{j} q \overline{l}} + R_{i
\overline{l} p \overline{q}} R_{q \overline{p} k \overline{j}} + R_{i \overline{j} k \overline{l}} \nonumber\\
& &  \;\;\; -{1\over 2} \left( R_{i \overline{p}}R_{p \overline{j}
k \overline{l}}  + R_{p \overline{j}}R_{i \overline{p} k
\overline{l}} + R_{k \overline{p}}R_{i \overline{j} p
\overline{l}} + R_{p \overline{l}}R_{i \overline{j} k
\overline{p}} \right). \label{eq:evolutio of curvature1}
\end{eqnarray}

The evolution equation for Ricci curvature and scalar curvature
are
 \begin{eqnarray} {{\p R_{i \b j}}\over {\p t}} & = & \triangle
  R_{i\b j} + R_{i\b j p \b q} R_{q \b p} -R_{i\b p} R_{p \b j},\label{eq:evolutio of curvature2}\\
  {{\p R}\over {\p t}} & = & \triangle R + R_{i\b j} R_{j\b i}- R.
  \label{eq:evolutio of curvature3}
  \end{eqnarray}

For
direct computations and using the evolved frames,we can obtain the
following evolution equations for the bisectional curvature:
 \begin{equation}
\pd {R_{i\bar jk\bar l}}{t} =\Delta R_{i\bar jk\bar l}- R_{i\bar
jk\bar l}+R_{i \bar j m\bar n}R_{n\bar m k\bar l}-R_{i\bar m k\bar
n}R_{m\bar j n\bar l}+R_{i\bar l m\bar n}R_{n\bar m k\bar l}
\label{eq:evolution of curvature4}
\end{equation}

As usual, the flow equation (\ref{eq:kahlerricciflow}) or
(\ref{eq:flowpotential}) is referred as the K\"ahler Ricci flow on
$X$. It is proved by Cao \cite{[2]}, who followed Yau's celebrated
work \cite{[19]}, that the K\"ahler Ricci flow exists globally for
any smooth initial K\"ahler metric.

\vskip8pt
\section{The traceless bisectional curvature operator}
\subsection{Definition and the evolution equations}
In Riemanian geometry ,  the curvature tensor for Riemannian
metric  can always be decomposed orthogonally into three parts:
$Rm=W+V+U,\;$  where $W$ is the Weyl tensor and $V,U$ are the
traceless Ricci part and the scalar curvature part respectively.
The decomposition for K\"ahler case is slight different.  The
bisectional curvature tensor can also be decomposed into
orthogonal parts as well.

Set
\begin{eqnarray}\twoten Sij&=&\twoten Rij-\frac 1n R \twoten gij = R_{i\b j}^0 ,\\
\four Sabcd &=&\four Rabcd -\frac 1n(\twoten Sab\twoten
gcd+\twoten Scd\twoten gab)-\frac 1{n^2}R\twoten gab\twoten gcd.
\label{eq:definition of S}
\end{eqnarray}

As in the Riemanian case,  the "Weyl" part $\four Sabcd $ is  also
trace free:$$\four Sabcd=\four Scdab,\;g^{a\bar b}\four Sabcd=0.$$
As in the previous subsection, under some evolved moving frame,
we can rewrite the evolution equation for curvature as below

\begin{prop}Along the K\"ahler Ricci flow the evolution
equation related the traceless bisectional curvature operator are
as follows:
 \begin{eqnarray}
\pd Rt&=&\Delta R-R+\frac 1nR^2+S_{\al\bar \bb}S_{\bb\bar \al}\\
\pd {\twoten Sab}{t}&=&\Delta \twoten Sab+\frac 1n (R-
n)\twoten Sab+S_{a\bar b i\bar j}S_{j\bar i}\\
\pd {\four Sabcd}{t}&=&\Delta S_{a\bar b c \bar d}- S_{a \bar b
c\bar d}+S_{a\bar b i\bar j}S_{j\bar i c\bar d}+S_{a\bar i j\bar
d}S_{i \bar b c\bar j} -S_{a\bar i c\bar j}S_{i\bar bj\bar
d}+\frac 1n\twoten Sab\twoten Scd \end{eqnarray} \end{prop}

  The bisectional curvature operator can be viewed as a symmetric operator
on the space of real $(1,1)$ forms $\Lambda^{1,1}(X).\;$   For any
pair of $(1,1)$ forms $\eta, \tau$, the action of the bisectional
curvature is:
\[
\cR(\eta, \tau) = R_{i\b j k \b l} \;\eta_{a \b b}\; \tau_{c\b
d}\; g^{i \b b} g^{a\b j} g^{k \b d} g^{c \b l}.
\]
If we decompose the space $\Lambda^{1,1}(X)$ into the line
consists of the multpile of the K\"ahler form and its orthogonal
completementary subspace $\Lambda^{1,1}_0(X),\;$  then the action
of $S_{i\b j k \b l}$ preserves $\Lambda^{1,1}_0(X).\;$  Denote
the action of $S_{i\b j k \b l}$ by $\cS.\;$  In some special basis, we will use $M$
to  denote the matrix of the
operator
$\cS.\;$We often referred
$\cS$ as the traceless bisectional curvature operator.  Moreover,
there is a nice decomposition formula for the bisectional
curvature operator in $ \Lambda^{1,1}(X):\;$
\begin{equation}
\left( \begin{array}{lcl}  R &Ric^0\\
{Ric^0}^t & \cS\end{array}\right).
\end{equation}

If the action of $\cS$ in $\Lambda^{1,1}_0(X),$ is non-negative,
then we call the underlying K\"ahler metric has a {\it
non-negative} traceless bisectional curvature operator.    If the
action of $\cS$ in $\Lambda^{1,1}_0(X)$ has a property that the
sum of any two eigenvalues is non-negative, then we call the
underlying K\"ahler metric has a {\it 2-non-negative} traceless
bisectional curvature operator.

%In this paper, we study the traceless bisectional curvature operator $\cS:\;$ both its
%algebraic properties and its behaviour under the K\"ahler Ricci flow.

\subsection{Geometric properties of the traceless  bisectional curvature operator}
In this subsection, we derive some geometric  properties of the traceless
bisectional curvature operator.  First,  in any local coordinate,
after fixing an frame such that the metric tensor at the origin is
identity matrix.   There is  a natural orthonormal  basis for
$\Lambda^{1,1}_0(X)$
 at the origin point (here $i, j =1, 2,\cdots n$):
$$\{\sqrt{-1}dz^i\wedge d\bar
z^i,dz^i\wedge d\bar z^j-dz^j\wedge d\bar z^i,\sqrt{-1}(dz^i\wedge
d\bar z^j+dz^j\wedge d\bar z^i)\}.$$ Note that the space of traceless
real (1,1) form is spanned by
$$\{{\sqrt{-1}}(dz^i\wedge d\bar z^i-dz^j\wedge d\bar
z^j),dz^i\wedge d\bar z^j-dz^j\wedge d\bar
z^i,\sqrt{-1}(dz^i\wedge d\bar z^j+dz^j\wedge d\bar z^i)\}.$$ For
convenience, we use the following spaces in this paper:

\begin{defi}:  The space $\La^{1,1}_0(X)$ is locally spanned by  the
following elements: \beqs A^{ij}&=&dz^i\wedge d\bar z^i-dz^j\wedge
d\bar z^j, \\B^{ij}&=&dz^i\wedge d\bar z^j+dz^j\wedge d\bar z^i,
\\C^{ij}&=&-\sqrt{-1}(dz^i\wedge d\bar z^j-dz^j\wedge d\bar z^i)
\eeqs and  the space $\La^{1,1}(X)$ is locally spanned by the
following elements: \beqs a^{ii}&=&2dz^i\wedge d\bar z^i, \\
B^{ij}&=&dz^i\wedge d\bar z^j+dz^j\wedge d\bar
z^i,\\C^{ij}&=&-\sqrt{-1}(dz^i\wedge d\bar z^j-dz^j\wedge d\bar
z^i), \eeqs
\end{defi}
where $i, j = 1,2 \cdots n.\;$  Our definition differs from the
space of the traceless (1,1) form because we want the
eigenvalues of $\cS$  to be positive for Fubni-Study metric.

\begin{prop}
If the traceless bisectional curvature operator is 2-nonnegative,
then the orthogonal  bisectional curvature is nonnegative.  If the
traceless bisectional curvature operator is nonnegative, then we
have the following inequalities:
$$
R_{i\bar i i\bar i}+R_{j\bar j j\bar j}\geq 2R_{i\bar i j\bar
j}\geq 0 ,\;\;R_{i\bar i}+R_{j\bar j}\geq 0
$$ for any\;\;$i\neq j$.
\end{prop}

\begin{proof}.(1) If $A$ is a symmetric matrix and the sum of two
lowest eigenvalues of $A$ is nonnegative,  then
$A_{ii}+A_{jj}\geq0$.  To see this, assume $m_1\leq m_2\leq\cdots
\leq m_n $ are the eigenvalues of $A$, then we have
$$
m_1+m_2=inf\{A(x,x)+A(y,y)| |x|=|y|=1,x\bot y \}\geq0 ,
$$
so we have
$$
A_{ii}+A_{jj}=A(e_i,e_i)+A(e_j,e_j)\geq0
$$
where $\{e_i\}$ are the standard basis of $\RR^n$.

(2)Since the matrix of $\cS$ is the same as the matrix of
curvature operator $Rm$ when acting on the space $\La_0^{1,1}(X)$,
so
$$
R(B^{ij},B^{ij})+R(C^{ij},C^{ij})\geq0,
$$
then simplify the above formula we have
$$
\four Riijj\geq0,\;\;\forall i\neq j.
$$

(3)If the traceless bisectional curvature operator is nonnegative,
Since the matrices of $\cS$ and $Rm$  are the same under the basis
of $\{A^{1i},B^{ij},C^{ij}\}$, then $R(A^{1i},A^{1i})\geq0$ and
$$
\left |\begin{array}{cc}
 R(A^{1i},A^{1i})&R(A^{1i},A^{1j})  \\
 R(A^{1i},A^{1j})&R(A^{1j},A^{1j})
\end{array}\right |\geq0.
$$
So
$$
R(A^{1i},A^{1i})+R(A^{1j},A^{1j})\geq
2\sqrt{R(A^{1i},A^{1i})R(A^{1j},A^{1j})}\geq2R(A^{1i},A^{1j}).
$$
Then we have
$$
R(A^{ij},A^{ij})=R(A^{1i},A^{1i})+R(A^{1j},A^{1j})-2R(A^{1i},A^{1j})\geq0.
$$
i.e.  $ R(A^{ij},A^{ij})=\four Riiii+\four Rjjjj-2\four Riijj\geq0
$.  Thus
\begin{eqnarray*}
R_{i\bar i}+R_{j\bar j}&=&\four Riiii+\sum_{\al\neq i} \four
R{\al}{\al}ii+\four Rjjjj+\sum_{\bb\neq j} \four R{\bb}{\bb}jj\\
&\geq&2\four Riijj+\sum_{\al\neq i} \four
R{\al}{\al}ii+\sum_{\bb\neq j} \four R{\bb}{\bb}jj
\\&\geq&0.
\end{eqnarray*}
 where $i\neq j$.  This finishes the proof of the proposition.
\end{proof}

\section{Proof of Theorem 1.2}

We follow notations in the previous section.  Note that
proposition 3.3 already show the first and last parts of Theorem
1.2. We need a technical lemma first.
\begin{lem}
If $M$ is a symmetric $m\times m$  matrix and satisfies \\
(1)$\sum_i M_{ij}=0,\;(\forall 1\leq j\leq m)$; \\
(2)$M_{ij}+M_{(i+1)(j+1)}-M_{i(j+1)}-M_{(i+1)j}=0, \;(\forall
1\leq i,j\leq m-1)$;\\then $M=0$ .
\end{lem}
\begin{proof}.From (1), we know
\begin{equation}\sum_i M_{ii}+2\sum_{i<j}M_{ij}=0.\label{eq:m1}
\end{equation} And from (2),we have $M_{ij}+M_{kl}-M_{il}-M_{jk}=0.\;$  So
$M_{ii}+M_{jj}-2M_{ij}=0 (\forall i<j).\;$  Then we have
\begin{equation}
(n-1)\sum_i M_{ii}-2\sum_{i<j}M_{ij}=0.\label{eq:m2}
\end{equation}
Thus from (\ref{eq:m1}) and (\ref{eq:m2}) we have
\begin{equation}
\sum_i M_{ii}=0,\;\sum_{i<j}M_{ij}=0.\label{eq:m3}
\end{equation}
From (\ref{eq:m3}) and $M_{ij}=M_{1i}+M_{1j}-M_{11}$,  we have
$\sum_{1<i<j} (M_{1i}+M_{1j}-M_{11})=0.\;$  Thus
$$(m-2)\sum_{i>1} M_{1i}-\frac {(m-1)(m-2)}{2}M_{11}=0.$$
Since $\sum_{i>1} M_{1i}+M_{11}=0$, we have $M_{11}=0$ if $m\geq
3$. Now $M_{ij}=M_{1i}+M_{1j}$, from (1) we have
$$0=\sum_j M_{ij}=\sum_j (M_{1i}+M_{1j})=mM_{1i}+\sum_j M_{1j}=mM_{1i}.$$
So $M_{1i}=0$, then $M_{ij}=0$.  We can easily verify that the
lemma still holds when $m=2$.  This finishes the proof of the
lemma.
\end{proof}

Now we are ready to give a proof of Theorem 1.2.
\begin{proof}.  We divide the proof into three parts.

 (1) In this part,  we want to prove that the trace of $\cS$ is bounded
  by the scalar curvature.  If the traceless bisectional curvature operator is
nonnegative, the proof is easy.  In this case,
Proposition 3.3 implies that

\begin{equation} R_{i\bar i i\bar i}+R_{j\bar j j\bar j}\geq 2R_{i\bar i
j\bar j}\geq 0,\;\forall i\neq j.\label{eq:Riijj}
\end{equation}
Let $M$ be the matrix expression of $\cS$ acting on $\La_0^{1,1}(X)$ with respect
to the basis (c.f. Definition 3.2). The
trace of $M$ is
\begin{eqnarray}
{\rm tr}(M)&=&\sum_{2\leq j\leq n}
R(A^{1j},A^{1j})+\sum_{i<j}(R(B^{ij},B^{ij})+R(C^{ij},C^{ij})) \nonumber\\
&=&\sum_{2\leq j\leq n} (\four R1111+\four Rjjjj-2\four
R11jj)+4\sum_{i<j}\four Riijj \nonumber\\
&=&R+(n-2)\four R1111+2\sum_{2\leq i<j}\four Riijj.
\label{eq:trace}
\end{eqnarray}
 From (\ref{eq:Riijj})(\ref{eq:trace}),  when $n\geq 4$
\beqs {\rm tr}(M)
&\leq &R+(n-2)\four R1111+\frac {n-2}{2}\sum_{2\leq i<j}2\four Riijj\\
&\leq &R+(n-2)\sum_{1\leq k\leq n} \four Rkkkk\\
&\leq &(n-1)R.\eeqs We can easily check that the above inequality
still holds when $n=2,3$.  In other words,  the trace of $\cS$ can
be bounded by the scalar
curvature.\\

If the traceless bisectional curvature operator is 2-nonnegative,
we need to bound $\four R1111$ and $\four Riijj$ from
(\ref{eq:trace}). Choose the basis $\{A^{1j},B^{ij},C^{ij}\}$ of
the space $\La_0^{1,1}(X) $ as in Subsection 3.2.  Then
$$R(A^{1i},A^{1i})+R(A^{1j},A^{1j})\geq 0,$$ i.e.,
$$\four R11ii+\four R11jj\leq \four R1111+\frac 12\four Riiii+\frac 12 \four Rjjjj.$$
Similarly, if we choose the basis $\{A^{kj},B^{ij},C^{ij}\},$  we
have
\begin{eqnarray}
\four Riijj+\four Riikk&\leq& \four Riiii+\frac 12\four
Rjjjj+\frac 12 \four Rkkkk,\nonumber\\\four Rjjii+\four
Rjjkk&\leq& \four Rjjjj+\frac 12\four Riiii+\frac 12 \four
Rkkkk,\nonumber\\\four Rkkii+\four Rkkjj&\leq&\four Rkkkk+\frac
12\four Rjjjj+\frac 12 \four Riiii.\label{eq:eqR}
\end{eqnarray}
Consequently,
$$\sum_{k\leq l}\four Rkkll\leq \sum_{k}\four Rkkkk .$$
Therefore,
$$
R=\sum_{k}\four Rkkkk+2\sum_{k<l}\four Rkkll\leq 3\sum_k \four
Rkkkk\leq 3R.
$$
The last inequality follows from $\four Rkkll\geq0  \;(\forall
k\neq l)$. Consequently,
$$\frac R3\leq \sum_k \four Rkkkk\leq R.$$
Since all the orthogonal bisectional curvature is nonnegative,  all of them are uniformly bounded by
 the scalar curvature.  Let us assume that the
holomorphic sectional curvature satisfies the following inequality
$$\four R1111\leq \four R2222\leq \cdots \leq \four Rnnnn.$$
We claim that $\four R3333\geq 0$.  To see this, if $\four R3333<0,$
then
$$\four R1111\leq \four R2222<0.$$
Consequently, we have $$M_{11}=\four R1111+\four R2222-2\four
R1122<0,$$
$$M_{22}=\four R1111+\four R3333-2\four R1133<0.$$
    Thus,  $M_{11}+M_{22}<0$,  which  is a contradiction! Consequently, we have
\[\four R3333\geq 0.\]
We want to again devide into three cases for discussions.

\textbf{Case 1}. If $\four R1111\geq 0$, then
$$0\leq \four R1111\leq \sum_k \four Rkkkk\leq R.$$

\textbf{Case 2}. If $\four R1111\leq \four R2222< 0$.  Since $M$
is 2-nonnegative and from (\ref{eq:eqR}), we have
 \beqs 2\four R1111+\four R2222+\four R3333&\geq& 0,\\
 2\four R2222+\four R1111+\four R3333&\geq &0.
 \eeqs
Thus
$$
|\four R1111|+|\four R2222|\leq \frac 23\four R3333.
$$
Therefore,
\begin{equation}
R\geq \four R3333-(|\four R1111|+|\four R2222|)\geq  \frac 13\four
R3333.\label{eq:formula1}
\end{equation}
Moreover, from $M_{11}+M_{22}\geq 0$, We have
\beqs M_{11}=\four R1111+\four
R2222-2\four R1122&<&0 ,\\M_{22}=\four R1111+\four R3333-2\four
R1133&>&0. \eeqs Consequently,
\begin{equation}
2\four R1111+\four R3333\geq 2\four R1111+\four R2222+\four
R3333\geq 2(\four R1122+\four R1133)\geq 0.\end{equation}

Thus, \begin{equation}|\four R1111|\leq \frac 12\four
R3333.\label{eq:formula2}\end{equation}

Combining inequalities  (\ref{eq:formula1}) and (\ref{eq:formula2}), we have
$$|\four R1111|\leq \frac R6.$$

\textbf{Case 3}.  If $\four R1111<0\leq \four R2222$, and if
$M_{11}=\four R1111+\four R2222-2\four R1122\geq 0,$  then
$$|\four R1111|\leq \four R2222.$$  From the proof of case 2, we
know $|\four R1111|\leq C(n)R$.  If $M_{11}<0$, then
$M_{22}\geq0$. So we have
$$|\four R1111|\leq \four R3333.$$Thus,  we can bound $|\four
R1111|$\;by the scalar curvature in this case.

In summary,  we can always bound ${\rm tr}(M) $ by the scalar
curvature $R$ if the bisectional curvature is 2-nonnegative.

(2) In this part,  we want to prove that  every $\four Sijkl$ can be
uniquely represented by the entries of the matrix $M$.  This is
equivalent to say that $\four Sijkl=0$ if $M=0$.  Assume $M=0$.\\
(Calculate $\four Siiii,\four Siijj$) Assume $ T_{ij}=\four
Siijj, T=(T_{ij}). $ We want to find $T$ from $M$.  Note that \beqs
T_{ij}+T_{kl} - T_{il}-T_{kj}&=&S(dz^i\wedge d\bar z^i-dz^k\wedge
d\bar z^k,dz^j\wedge d\bar z^j-dz^l\wedge d\bar
z^l)\\&=&S(A^{1k}-A^{1i},A^{1l}-A^{1j})\\&=&0\eeqs since
$S(A^{1k},A^{1l})=0$.   Thus, the following equations hold:
\begin{eqnarray*}
T_{ij}+T_{(i+1)(j+1)}-T_{i(j+1)}-T_{(i+1)j}&=& 0\\
\sum_i T_{ij}&=&0,  \;\forall j  =1,2 \cdots n
\end{eqnarray*}
From Lemma 4.1,  we have $T=0$, i.e. \begeq \four Siiii=\four
Siijj=0. \label{eq:S11}\endeq (Calculate other $\four Sijkl$)
Since $dz^i\wedge d\bar z^j=\frac 12(B^{ij}+\sqrt{-1}C^{ij}),$ we
have \begeq \four Sijij=0(\forall i\neq j).\label{eq:S2}\endeq
From $ S(B^{ij},B^{ik})=0,S(C^{ij},C^{ik})=0$ and $
S(B^{ij},C^{ik})=0,S(C^{ij},B^{ik})=0, $  where $i\neq j,j\neq
k,k\neq i,$  we have $\four Sijik +\four Sjiki=0,\four Sijik
=\four Sjiki. $ Consequently,   \begeq \four Sijik=\four
Sjiki=0.\;\label{eq:S3}\endeq From
$S(A^{ij},B^{kl})=0,S(A^{ij},C^{kl})=0, $ where $i\neq j,k\neq l,$
we have $\four Siikl=\four Sjjkl.\;$  Since $\sum_i \four
Siikl=0,$ we have \begeq \four Siikl=0.\label{eq:S4}\endeq  Let
$l=i$, we have \begeq \four  Siiki=0.\label{eq:S5}\endeq From
(\ref{eq:S3})and $S(A^{ij},B^{ij})+\sqrt{-1}S(A^{ij},C^{ij})=0$
where $i\neq j,$ we have \begeq \four Siiij=\four
Sjjij=0.\label{eq:S6}\endeq
 From
$S(B^{ij},B^{ij})=0$ where $i\neq j$, we have \begeq \four
Sijji=0.\label{eq:S7}\endeq From (\ref{eq:S2}) and
$S(B^{ij},B^{kl})=0$ where $i\neq j,j\neq k,k\neq
l,l\neq i$, we have \begeq \four Sijkl=0.\label{eq:S8}\endeq\\
From (\ref{eq:S11})-(\ref{eq:S8}), we know
$$\four Sijkl=0.$$

(3) Finally we prove that the bisectional curvature is bounded by the
scalar curvature.  Let $m_1\leq m_2\leq \cdots \leq m_s (s=n^2-1)$
be the eigenvalues of the matrix $M$.  If $M$ is positive, then
all the eigenvalues can be bounded by the scalar curvature.
Otherwise since  $M$ is 2-positive, the eigenvalues must satisfy
the following inequalities:
$$m_1\leq 0\leq m_2\leq \cdots \leq m_s,\;m_1+m_2>0.$$
Then,  $m_3< \mathit{tr}(M)$ and $|m_1|<m_2\leq m_3<\mathit{tr}(M).\;$
In other words,  all the eigenvalues can be bounded by the scalar
curvature. \\

Since all the eigenvalues of the matrix $M$ are bounded by $R$,
 all of its entries of $M$ are also  bounded by the scalar curvature function.
 From (2), every $\four Sijkl$ is bounded by $R$, i.e.  $|\four Sijkl|<C(n)R$.
Thus $|S_{k\bar l}|\leq C(n)R$.  From the definition of $\cS$
(c.f. (\ref{eq:definition of S}), we have $$|\four Rijkl|\leq
C(n)R.$$

\end{proof}

\section{Proof of Theorem 1.1}
 In this section ,we are ready to  prove  Theorem
 1.1.  Note that in \cite{[9]}, the
  positive curvature operator is preserved  and in \cite{[4]} the 2-positive curvature operator is
  preserved along the Ricci  flow, one can  also see both from  \cite{[10]}.  Our proof
  here is similar to theirs. First recall R. Hamilton's lemma:
\begin{lem} {\rm {\cite{[9]}}}. Suppose $\partial f/\partial t=\Delta f+\phi(f).$
Let $s(f)$ be a concave function on the bundle invariant under
parallel translation whose level curves $s(f)\leq c$ are preserved
by the ODE $df/dt=\phi(f)$.  Then the inequality $s(f)\leq c$ is
preserved by the PDE for any constant $c$.  Furthermore if
$s(f)<c$ at one point at time $t=0$, then $s(f)<c$ everywhere on
$X$ for all $t>0$.

\end{lem}

Now we begin to prove Theorem 1.1.
%\subsection{Proof of Theorem 1.1 }
\begin{proof}.
(1) Define
%\phi^{\rho}_{a\bar b}
$$[\phi^\la,\phi^\mu]_{a\bar b}=\phi^{\la}_{a\bar m}\phi^{\mu}_{m\bar b}-\phi^{\mu}_{a\bar m}\phi^{\la}_{m\bar b}=C^{\la\mu}_{\rho}\phi^{\rho}_{a\bar b}.$$
By calculation we have the following relations: \beqs
[A^{ij},B^{ij}]&=&2\sqrt{-1}C^{ij},[A^{ij},B^{ik}]=\sqrt{-1}C^{ik},[A^{ij},B^{jk}]=\sqrt{-1}C^{kj},\\
{[}B^{ij},C^{ij}]&=&2\sqrt{-1}A^{ij},[B^{ij},C^{ik}]=-\sqrt{-1}B^{jk},[B^{ij},C^{jk}]=-\sqrt{-1}B^{ik},\\
{[}C^{ij},A^{ij}]&=&2\sqrt{-1}B^{ij},[C^{ij},A^{ik}]=\sqrt{-1}B^{ji},[C^{ij},A^{jk}]=-\sqrt{-1}B^{ij},\\
{[}B^{ij},B^{ik}]&=&\sqrt{-1}C^{jk},[C^{ij},C^{ik}]=\sqrt{-1}C^{jk}.
\eeqs where $i\neq j,j\neq k,k\neq i$,  and other Lie brackets are
zero.  Note that all $C^{\la\mu}_{\rho}$ are zeros or pure
imaginary numbers. \beqs \four
Samnd\four Smbcn&-&\four Samcn\four Smbnd\\
&=&M_{\al\bb}\phi^{\al}_{a\bar m}\phi^{\bb}_{n\bar
d}M_{\gamma\delta}\phi^{\gamma}_{m\bar b}\phi^{\delta}_{c\bar
n}-M_{\al\bb}\phi^{\al}_{a\bar m}\phi^{\bb}_{c\bar
n}M_{\gamma\delta}\phi^{\gamma}_{m\bar b}\phi^{\delta}_{n\bar d}\\
&=&M_{\al\bb}M_{\gamma\delta}\phi^{\al}_{a\bar
m}\phi^{\gamma}_{m\bar b}(\phi^{\bb}_{n\bar d}\phi^{\delta}_{c\bar
n}-\phi^{\bb}_{c\bar n}\phi^{\delta}_{n\bar d}) \\
&=&M_{\al\bb}M_{\gamma\delta}\phi^{\al}_{a\bar
m}\phi^{\gamma}_{m\bar b}C^{\delta\bb}_{\rho}\phi^{\rho}_{c\bar
d}\\
&=&-\frac
12C^{\al\gamma}_{q}C^{\bb\delta}_{p}M_{\al\bb}M_{\gamma\delta}\phi^{q}_{a\bar
b}\phi^{p}_{c\bar d}.\eeqs Define
$$M^{\#}_{qp}=C^{\al\gamma}_{q}C^{\bb\delta}_{p}M_{\al\bb}M_{\gamma\delta},\eqno (4.1)$$
then we have
$$\four
Samnd\four Smbcn-\four Samcn\four Smbnd=-\frac
12M^\#_{qp}\phi^{q}_{a\bar b}\phi^{p}_{c\bar d}$$ and $$\pd
{M}{t}=-M+M^2-\frac 12M^\#+\frac 1nT .\eqno (4.2)$$ \vskip0.5cm
Now we have the following lemma:

\begin{lem} If all $C^{\al\gamma}_{q}$  are  real and
$M\geq0$, then $M^{\#}\geq 0$.
\end{lem}

\begin{proof}.
 Without loss of generality, we may choose a basis $\{\phi^{\al}\}$ which
diagonalizes $M$, so that $M_{\al\bb}=\delta_{\al\bb}M_{\al\al}$.
For any $v=v^{\al}\phi^\al$, we have $$M^{\#}(v,v)=(v^\al
C^{ai}_{\al})(v^\bb C^{bj}_{\bb})M_{ab}M_{ij}=(v^\al
C^{ai}_{\al})^2M_{aa}M_{ii}\geq 0.$$ The lemma is then proved.\end{proof}

Now we return to the proof of Theorem 1.1 again.  Since in our
case all $C^{\la\mu}_{\rho}$ are zero or pure imaginary numbers,
then $M^\#\leq0$  if $M\geq0$. Since $T$ is always non-negative, we have
$$\pd {M}{t}=-M+M^2-\frac 12M^\#+\frac 1nT \geq 0.$$ when $M=0$.
Note that $M\geq 0$ is convex and  $M(0)\geq0$, we have
$M(t)\geq0$ for all $t>0$.  In other words, the nonnegative
traceless bisectional curvature operator is preserved.  From Lemma
5.1, if $M$ is positive at one point at time $t=0$, then $M$ is
positive everywhere for all time $t>0$.
\\

(2) We want to prove that the 2-nonnegative traceless bisectional curvature
operator is preserved along the K\"ahler Ricci flow.  Let us  assume that
the eigenvalues of the traceless bisectional curvature operator on
$\La_0^{1,1}(X)$ are $\la_1\leq \la_2\leq \cdots\la_m$, where
$m=n^2-1$. From (4.1)(4.2), we have

$$\beqa{ll}
 \frac {d}{dt}(\la_1+\la_2)&\geq \frac
{d}{dt}(M_{11}+M_{22})\\
&\geq-(\la_1+\la_2)+(\la_1^2+\la_2^2)-\frac
12\sum_{p,q}((C_1^{pq})^2+(C_2^{pq})^2)\la_p\la_q.
\eeqa\eqno(4.3).$$

Note that the right side \beqs\frac
12\sum_{p,q}((C_1^{pq})^2+(C_2^{pq})^2)\la_p\la_q
&=&\sum_{p<q}((C_1^{pq})^2+(C_2^{pq})^2)\la_p\la_q\\&
=&\sum_{q\geq 3}(C_1^{2q})^2 (\la_1+\la_2)\la_q+\sum_{p,q\geq
3}((C_1^{pq})^2+(C_2^{pq})^2)\la_p\la_q .\eeqs

Note that $\la_m\geq \cdots\geq \la_2\geq 0.\;$  If
$\la_1+\la_2=0$, then the right side of (4.3) is nonnegative.
Since $\la_1+\la_2$ is a concave function on $X$,
 $\la_1+\la_2\geq 0$ is preserved.   From Lemma
5.1, if $\la_1+\la_2>0$ is positive at one point at time $t=0$,
then $\la_1+\la_2>0$ is positive everywhere for all time
$t>0$.\end{proof}

(3)  Now we prove the last part of Theorem 1.1. If the traceless
bisectional curvature positive is non-negative or 2-non-negative,
Theorem 1.2, implies that
$$\four Riijj\geq0 ,\;\forall i\neq j.$$ Let us assume initially the Ricci curvature is positive and
after finite time $t_0> 0$, at some point $p\in X$
$R_{i\bar j}$ vanishes at least at one direction.  For
convenience, set this direction as $\frac {\partial}{\partial
z_1}$ and diagonalize the Ricci curvature at this point.  Then
$$\pd {R_{1\bar 1}}{t} |_{t_0}\geq \four R11jj R_{j\bar j}-R_{1\bar 1}R_{1\bar 1}=\sum_{j=2}^n \four R11jj R_{j\bar j}\geq  0.$$
By Hamilton's maximum principle for tensors, this is enough to show that the positivity
of  Ricci curvature is preserved under the condition.

\section{Proof of Corollary 1.4 }
We only need
  to prove the K\"ahler Ricci flow convergence by sequences to some K\"ahler Ricci
  soliton when the bisectional curvature is uniformly bounded from
  Theorem 1.2.  In \cite{[17]},  She proved that $\tau$-flow converges by sequence to some Ricci
  soliton when the  curvature operator and the diameter are  uniformly bounded.  In
  \cite{[18]},   she proved that the K\"ahler Ricci flow converges
by sequence to some K\"ahler Ricci  soliton except a set of
isolated points on any complex compact K\"ahler surface without
any curvature assumptions.  First let us recall Perelman's no
local collapsing theorem:
\begin{defi}{\rm \cite{[13]}}.  Let $g_{ij}(t)$ be a smooth solution to the Ricci flow
$(g_{ij})_t$ on $[0,T)$ on a Riemannian manifold  X of dimension
$n$. We say that $g_{ij}(t)$ is locally collapsing at $T$, if
there is a sequence of times $t_k\rightarrow T$ and a sequence of
metric balls $B_k=B_k(p_k,r_k)$ at times $t_k$, such that $\frac
{r_k^2}{t_k}$ is bounded,  $|Rm|(g_{ij})(t_k)\leq r_k^{-2}$ in
$B_k$ and $\frac {\mathit{Vol}(B_k)}{r_k^{n}}\rightarrow 0$.
\end{defi}

\begin{lem}{\rm \cite{[13]}}.  If $X$ is closed and $T<\infty$, then
$g_{ij}(t)$ is not locally collapsing at $T$.
\end{lem}

Now we begin to prove Corollary 1.4.
\begin{proof}.  We divide the proof into two parts.\\
(1) First we are ready to prove that the injectivity radii have a
uniformly  positive lower bound along the K\"ahler Ricci flow. If
the traceless bisectional curvature operator is 2-nonnegative and
the scalar curvature is bounded along the flow, Theorem 1.2
implies that the curvature tensor is uniformly bounded.

% Consequently,  the diameter is uniformly bounded along the
%K\"ahler Ricci flow.   Since the volume is constant along the
%K\"ahler Ricci flow,  we have uniformly positive lower bound on
%the injectivity radii along the flow.  \\

\begin{claim}The injectivity radius has a uniformly positive lower
bound along the flow.
\end{claim}
{\it Proof.} Let $(X,g_{i\bar j})$ be the  K\"ahler Ricci flow.
Fix $T>0$.  Now we re-scale the metric
 \begin{equation} \bar g_{i\bar
j}(s)=(T-s)g_{i\bar j}(-\log (\frac {T-s}{T})), \qquad s\in
[0,T).\label{eq:eqg}\end{equation}

Then, $\bar g_{i\bar j}(s)$ is a solution with finite maximal
existence interval to the unnormalized K\"ahler Ricci flow $\pd
{g_{i\bar j}}{s}=-R_{i\bar j}.$   Lemma 6.2 implies that $(X,\bar
g_{i\bar j}(s))$ is not locally collapsing. In other words, for
any sequence of times $s_k\rightarrow T,$  any  sequence of metric
balls $B_k=B_k(x_k,r_k)$ at times $s_k$, such that $\frac
{r_k^2}{s_k}$ is bounded and $|Rm|(\bar g_{ij})(s_k)\leq r_k^{-2}$
in $B_k$, there exists a constant $\delta>0$ such that

\begin{equation} \frac {\mathit{Vol}(B_k)}{r_k^{2n}}\geq
\delta.\label{eq:vol1}\end{equation}

Since $ |Rm|(g_{i\bar j}(t))$ is uniformly bounded along the
K\"ahler Ricci flow, for the un-normalized flow, we have $$|Rm|(\bar
g_{i\bar j}(s))\leq \frac C{T-s}.$$ We claim that there exists a
constant $\ee>0$ such that $\mathit{inj}(\bar g(s))\geq \sqrt{T-s}\;
\ee$. We prove this by contradiction.  If there exist a sequence of
times $s_k\rightarrow T$, such that
$$\frac {\mathit{inj}(\bar
g(s_k))}{\sqrt{T-s_k}}\rightarrow 0.$$

We re-scale the metric
$$h(s_k)=\frac 1{T-s_k}\bar g(s_k).$$
Let $r_k^2=T-s_k$, then  \begin{equation}  |Rm|(h(s_k))\leq
C,\mathit{inj}(h(s_k))\rightarrow 0.\label{eq:vol2}\end{equation}
From (\ref{eq:vol1}), we have \begin{equation}
\mathit{Vol}(B_{h(s_k)}(x_k,1))\geq \delta.\label{eq:vol3}
\end{equation} Then, (\ref{eq:vol2})(\ref{eq:vol3}) contradict with
J.Cheeger's injectivity radius estimate (c.f.\cite{{[14a]}}). Thus,
we have $\mathit{inj}(\bar g(s))\geq \sqrt{T-s} \;\ee$. Together
with (\ref{eq:eqg}), we have
$$\mathit{inj}(g(t))\geq \ee > 0.$$

\begin{claim}The diameter has a uniformly upper bound along the
flow.
\end{claim}
{\it Proof.} To see this,  we assume that there are $N$ points
$p_1,p_2,...p_N$ such that
$$dist_{g(t)}(p_i,p_j)\geq 2\ee,\qquad
\forall 1\leq i\neq j\leq N$$ where $\ee>0 $ is the uniformly lower
bound on the injectivity radius from Claim 6.3. Hence, the balls
$B_{g(t)}(p_i,\ee)$ are embedded and pairwisely disjoint. Since the
curvature operator is uniformly bounded and from the volume
comparison theorem,
$$V\geq \sum_{i=1}^N \mathit{Vol}(B_{g(t)}(p_i,\ee))\geq N C\ee^{2n}.$$Since the volume $V$ is
fixed along the flow,  $N$ is bounded  from
above.  Consequently the diameter  has a uniformly upper bound along the flow.\\

(2) Now we return to the proof of Corollary 1.4. Since we have
uniformly bounds on curvature tensor and uniformly lower bound on
the injectivity radius,  by Hamilton's compactness theorem, for
every $t_k\rightarrow \infty$ as $k\rightarrow \infty$, there exists
a subsequence such that $(X,g(t_k+t))$ converges to $(X,h(t))$, in
the sense that there exist diffeomorphisms $\phi_i:X\rightarrow X$,
such that $\phi_i^*g(t_k+t)$ converge uniformly together with their
covariant derivatives to metrics $h(t)$ on any compact subsets. For
every sequence of times $t_k\rightarrow \infty$, there exists a
subsequence, such that the $(X,g(t_k+t))$ converges to a K\"ahler
Ricci soliton as $k\rightarrow\infty$.  Now we outline that
process (c.f. \cite{[17]}\cite{[18]}).\\

Corresponding to the K\"ahler Ricci flow, Perelman's functional:
$$W(g,f,\frac 12)=(2\pi)^{-n}\int_X \;(|\Na f|^2+R+f-2n)e^{-f}dV_g.$$
One can show that $\mu(g,\frac 12)= \inf \{W(g,f,\frac
12)|f\;\mathit{satisfies}  \int_X \; (2\pi)^{-n}e^{-f}=1\}$ can be
achieved by a smooth function $f(t)$ such that $\mu(g,\frac
12)=W(g,f,\frac 12).\;$ Moreover,
$$2\Delta f-|\Na f|^2+R+f-2n=\mu(g,\frac
12).$$Then, we have the following results:\\
(i) $\lim_{i\rightarrow \infty}\; \mu(g_i,\frac 12)=\mu (h,\frac
12)$, so that $\lim_{i\rightarrow \infty}\;
\mu(g_i,\frac 12) $ is finite;\\
(ii) If $\mu(g_k,\frac 12)=W(g,f_k,\frac 12)$ and  $\mu(h,\frac
12)=W(h,f,\frac 12)$, then $f_k\rightarrow f$
 in $C^{2,\al}$ norm;\\
(iii) For any t, we can find $f_t$ such that $\mu(g(t),\frac
12)=W(g(t),f_t,\frac 12)$.  If we flow $f_t$ backward, we will get
a function $f_t(s)$.  Let $u_t=e^{-f_t}.$ Then,  $u_t(s)$
satisfies the following equation

$$\left \{\begin{array}{rl} \pd {}{s}u_t(s)&=-\Delta u_t(s)+(n-R)u_t(s)\\ u_t(t)&=e^{-f_t} \end{array} .\right.$$

Then, for every $A>0$, there exist a uniformly constant $C$
depending on $A$, such that $|u_t(s)|_{C^{2,\al}}\leq C$  for
every $t\geq A,s \in [t-A,t]$.\\

 Now in our case, fix $A>0$ and
assume $g_{i\bar j}-R_{i\bar j}=\partial_i\partial_{\bar j}u$,
\beqs \mu(g(t_i+A),\frac 12)-\mu(g(t_i),\frac 12)&\geq&
W(g(t_i+A),f_{t_i+A},\frac 12)-W(g(t_i),f_{t_i+A}(t_i),\frac
12)\\&=&(2\pi)^{-n}\int ^A_0\;|\Na\bar {\Na}
(f_{t_i+A}-u)|^2(t_i+s)dV_{g(t_i+s)} .\eeqs From (i), we have
$\lim_{i\rightarrow \infty}|\Na\bar {\Na}(f_{t_i+A}-u)|(t_i+s)=0$
for all most $x\in X$  and $s\in [0,A]$.  Since curvature is
bounded, $u(t_i+s)\rightarrow \bar u(s) $ in $C^{2,\al}$ norm for
any $s\in [0,A]$.  From  (iii) we can assume
$u_{t_i+A}(t_i+s)\rightarrow \bar u_1(s)$.  So we have $$\Na\bar
{\Na}(\bar f-\bar u)=0,$$where $\bar u_1=e^{-\bar f}$.  Then,
$\bar u(s)=\bar f(s)$ since both of them satisfy the same integral
condition
$$\int_X\;e^{-\bar f}dV_s=\int_X\;e^{-\bar u}dV_s=(2\pi)^n.$$
 Since $\bar f,\bar u$  satisfy the
following equations \beqs \pd{}{s}\bar f(s)&=&-\Delta \bar f+|\Na
\bar f|^2-R+n,
\\\pd{}{s}\bar u(s)&=&\Delta \bar u+\bar u+a(s),\\\Delta \bar u&=&n-R
 ,\eeqs we have $$\Delta \bar u-|\Na \bar u|^2+\bar
 u=-a(s).$$Thus $a(s)=-(2\pi)^{-n}\int_X\;\bar ue^{-\bar u}dV_s$.  Then, we can show that $a(s)$ is a constant.
 Note that $$\mu(h(t),\frac 12)=(2\pi)^{-n}\int_X \;(|\Na \bar u|^2-\Delta \bar u+\bar u-n)e^{-\bar u}dV_{g(t)}=-a(t)-n.$$
Since there exists a finite $\lim_{t\rightarrow\infty}\mu
(g(t),\frac 12)$, then $\mu(h(t),1/2)=\mu(h(s),1/2)$ for all
$s,t\in [0,A]$. Therefore,$a(t)$ is a constant.  Note that
$$\pd {}{t}\int_X\;(\bar u+a)dV_{h(t)}=-\int_X \;|\Na\Na u|^2dV_{h(t)}$$
and $$\pd {}{t}\int_X\;\bar udV_{h(t)}=\int_X\;(|\Na \bar
u|^2+\bar u\Delta \bar u)dV_{h(t)}=0.$$ Consequently  we have
$\Na\Na \bar u=0$.  Then, $h(t)$ is a K\"ahler Ricci soliton for
$t\in [0,A]$.
 Taking an monotone increasing sequence $A_j\rightarrow\infty$, we can
 obtain
functions $\bar u_{A_j}(t)$ for $t\in [0,A_j]$.  Then, we can show
that $\bar u_{A_j}=\bar u_{A_k}$ in $[0,A_j]$ for $j<k$.  Thus we
get a soliton for all time $t\geq 0$.
\end{proof}

Xiuxiong Chen,  xiu@math.wisc.edu, University of Wisconsin at
Madison, USA// Haozhao Li, lnsinx@math.pku.edu.cn,  Peking
University, Beijing, P. R. China

{\small\begin{thebibliography}{2}

\bibitem{[1]} S. Bando. On the three dimensional compact K\"ahler
manifolds of nonnegative bisectional curvature. \emph{J. D. G.},
19:283-297, 1984.

\bibitem{[2]} H. D. Cao. Deformation of K\"ahler metrics to
K\"ahler-Einstein metrics on compact K\"ahler manifolds.
\emph{Invent. Math.}, 81:359-372, 1985.
\bibitem{[2a]}
B. Chow; P. Lu. The maximum principle for systems of parabolic
equations subject to an avoidance set.  \emph{Pacific J. Math.}
214(2):201-222, 2004.

\bibitem{[3]} H.D. Cao, B. L. Chen  and X. P. Zhu. Ricci flow on K\"ahler
manifold of positive bisectional curvature, 2003. math.DG/0302087.

\bibitem{[4]} H. W. Chen.  Pointwise 1/4 -pinched 4-manifolds. \emph{Ann. Global
Anal. Geom.} 9(2):161-176, 1991.

\bibitem{[5]} X. X. Chen.  On K\"ahler
manifold with positive orthogonal bisectional curvature. In
prepation.

\bibitem{[6]} X. X. Chen,  G. Tian. Ricci flow on K\"ahler-Einstein
surfaces. \emph{Invent. Math.} 147(3):487-544, 2002.

\bibitem{[7]} X. X. Chen, G. Tian. Ricci flow on K\"ahler-Einstein manifolds.
arXiv:math. DG/0108179.

\bibitem{[7a]} D. Knopf.  Positivity of Ricci curvature under the K\"ahler-Ricci flow
arXiv:math. DG/0501108.

\bibitem{[8]} R. Hamilton. Three-manifolds with positive Ricci
curvature. \emph{J. Diff. Geom.}, 17:255-306, 1982.

\bibitem{[9]} R. Hamilton, Four manifolds with positive curvature
operator, \emph{J.Differential Geom.} 24:153-179, 1986.

\bibitem{[10]} R. Hamilton,The formation of singularities in the
Ricci flow,in \emph{Surverys in Differential Geom.}2, pp.7-136,
International Press, 1995.

\bibitem{[11a]} R. Hamilton,  Four-manifolds with positive isotropic curvature.
 \emph{Comm. Anal. Geom.} 5(1):1-92, 1997.

\bibitem{[11]} Ni, Lei. Ricci flow and nonnegativity of curvature.
arXiv:math.DG/0305246.

\bibitem{[12]} N. Mok. The uniformization theorem for compact K\"ahler
manifolds of nonnegative holomorphic bisectional curvature.
\emph{J. Differential Geom.}, 27:179-214, 1988.

\bibitem{[13]} G. Perelman.  The entropy formula for the Ricci flow and its
geometric applications, preprint arXiv:math.DG/0211159.

\bibitem{[14]} G. Perelman. K\"ahler-Ricci flow (unpublished work).

\bibitem{[14a]}Peter Petersen, Riemannian geometry,  Graduate texts
in mathematics, 171, Springer-Verlag, 1997.

\bibitem{[15]}Phong, Duong H., Sturm, Jacob. On the K\"ahler-Ricci flow on complex surfaces.
arXiv:math.DG/0407232.

\bibitem{[16]}Phong, Duong H.; Sturm, Jacob. On stability and the
convergence of the K\"ahler-Ricci flow. arXiv:math.DG/0412185.

\bibitem{[17]}N. Sesum.  Limiting behaviour of the Ricci flow, preprint
arXiv:math.DG/0402194.

\bibitem{[18]}N. Sesum. Convergence of a K\"ahler Ricci
flow;preprint.arXiv:math.DG/0402238.
\bibitem{[19]} S. T. Yau. On the Ricci curvature of a compact K\"ahler
manifold and the complex Monge-Ampere equation, I*, \it{Comm. Pure
Appl. Math.}, 31:339-441, 1978.
\end{thebibliography}}
\end{document}